\documentstyle[art10]{article}
\title{Canonical Geometrically Ruled Surfaces.}
\author{Luis Fuentes Garc\'{\i}a\thanks{Supported by an F.P.U.
fellowship of Spanish Government}
 \and{Manuel Pedreira P\'erez}
 }
\date{}

\newtheorem{teo}{Theorem}[section]
\newtheorem{defin}[teo]{Definition}
\newtheorem{prop}[teo]{Proposition}
\newtheorem{cor}[teo]{Corollary}
\newtheorem{lemma}[teo]{Lemma}

\newtheorem{rem}[teo]{Remark}

\newtheorem{propdef}[teo]{Proposition-Definition}
\font\euf=eufm10 at 12pt

\def\g2{\pi}
\def\L{{\cal L}}

\def\e{\mbox{\euf e}}
\def\eh{E}
\def\b{\mbox{\euf b}}
\def\bh{B}
\def\bb{\mbox{{\euf b}}}

\def\c{\mbox{\euf c}}
\def\ch{C}
\def\a{\mbox{\euf a}}
\def\ah{A}
\def\R{{\cal B}}
\def\aa{\mbox{\euf a}}
\def\k{{\cal K}}
\def\K{{\cal K}}
\def\gr{\deg}
\def\P{{\bf P}}

\def\Ram{{\bf R}}
\newcommand\p{{\bf P}}
\newcommand\E{{\cal E}}
\newcommand\Te{{\cal O}}
\newcommand\C{{C}}

\newcommand\A{\alpha}
\newcommand\B{\beta}

\def\qed{\hspace{\fill}$\rule{2mm}{2mm}$}
\newcommand\lrw{\longrightarrow}

\newcommand{\Div}{\mathop{\rm Div}\nolimits}

\newcommand{\Pic}{\mathop{\rm Pic}\nolimits}

\newcommand{\id}{\mathop{\rm id}\nolimits}

\newcommand{\Sym}{\mathop{\rm Sym}\nolimits}
\newcommand{\coker}{\mathop{\rm coker}\nolimits}
\newcommand{\Cliff}{\mathop{\rm Cliff}\nolimits}

\def\ov{\overline}
\def\vhi{\varphi}

\begin{document}
\maketitle

{\footnotesize{\bf Authors' address:} Departamento de Algebra, Universidad de Santiago
de Compostela. $15706$ Santiago de Compostela. Galicia. Spain. e-mail: {\tt
pedreira@zmat.usc.es}; \\ {\tt luisfg@usc.es}\\
{\bf Abstract:}  We prove the existence of canonical scrolls; that is, scrolls playing the
role of canonical curves. First of all, they provide the geometrical version of
Riemann Roch Theorem: any special scroll is the projection of a canonical scroll
and they allow to understand the classification of special scrolls in $P^N$.
Canonical scrolls correspond to the projective model of canonical geometrically
ruled surfaces over a smooth curve. We also prove that the
generic canonical scroll is projectively normal except in the hyperelliptic
case and for very particular cases in the nonhyperelliptic situation.\\
{\bf Mathematics Subject Classifications (1991):}
Primary, 14J26; secondary, 14H25, 14H45.\\ {\bf Key Words:} Ruled Surfaces, elementary
transformation, projective normality.}

\vspace{0.1cm}

{\Large\bf Introduction.}  A {\it geometrically ruled surface}, or simply a {\it
ruled surface}, is a $\P^1$-bundle over a smooth curve $X$ of genus $g>0$. It will be
denoted by $\pi: S=\P(\E_0)\lrw X$, where $\E_0$ is the associated normalized locally free sheaf of rank
$2$.  Let
$X_0$ be the minimum self-intersection curve of the ruled surface. Let $H\sim X_0+\b f$, $\b\in
\Div(X)$ be an unisecant divisor such that the complete linear system $|H|$ provides a regular morphism
$\phi_H:\P(\E_0)\lrw R\subset \P^N, n=\dim(|H|)$. We call $R$ a linearly normal scroll. We refer to
\cite{fuentes} for a systematic development of the projective theory of scrolls and ruled surfaces.

We define the {\it speciality} of a scroll $R$  as the superabun\-dance of $|H|$; that is,
$i(R)=s(\Te_S(H))=h^1(\Te_S(H))$. The scroll $R$ is called special if $i(R)>0$. If $W\subset \P^n$ is a
subspace and $\pi_W:R-(W\cap R)\lrw R'\subset \P^{n'}$ is the projection of $R$ from $W\cap R$, it is
well known that $i(R')-i(R)=\gr(W\cap R)-(dim(\langle W\cap R\rangle)+1)$, (see \cite{fuentes},
$\S5$), and the speciality of $R$ grows exactly the number of unassigned base points in the
linear system of hyperplane sections containing $W\cap R$. 

In the case of curves, the Riemann-Roch Theorem gives a nice geometrically interpretation of the
speciality of a nonhyperelliptic curve. Any linearly normal special curve $C$ is the projection  of a
canonical curve $C_{\K}$ from a set of points $A$. The canonical curve $C_{\K}$ is projectively normal
and it has speciality $1$. The speciality of $C$ is $\gr(A)-dim(\langle A \rangle)$.

In this paper we prove the existence of canonical scrolls; that is, the existence of scrolls
playing the role of canonical curves and, in particular, providing the geometrical version of
Riemann-Roch Theorem for ruled surfaces. Moreover, we prove that in general they are projectively normal
when the base curve is not hyperelliptic. 

Note that the study of special scrolls is equivalent to the study of special locally free sheaves of rank
$2$ over a smooth curve. The projection of a scroll corresponds to the elementary transformation of a
locally free sheaf of rank $2$ (see \cite{fuentes}, $\S4$ and \cite{maruyama}, $\S1$). In this way we
will prove the existence of locally free sheaves of rank $2$  and speciality $1$ such that any special
locally free sheaf of rank $2$ is obtained from them by elementary transformations.

The start point is a nice result mentioned by C. Segre in
\cite{segre1}: {\it if $R\subset \P^n$ is a linearly normal scroll and $\C\subset R$ is a bisecant curve
which has not double points out of the singular locus of $R$, then $\C$ is linearly normal and the
speciality of $\C$ is equal to the speciality of $R$}. We call $C$ a proper bisecant curve.

This theorem suggests the definition of canonical scroll: Let $X$ be a smooth curve of genus $g>0$ and
let $\C$ be a smooth curve of genus $\pi$ such that there exists an involution $\gamma:\C\lrw X$; that
is, a finite morphism of degree $2$. Let us suppose that $\C$ is
not hyperelliptic and it has genus $\pi\geq 3$. Let $\vhi_{\K}:\C\lrw \C_{\K}\subset
\P^{\pi-1}$ the canonical map. Then $\C_{\K}$ is a bisecant curve in the scroll
$S=\bigcup_{x\in X}\langle \gamma^{-1}(x)\rangle$. By Segre's result, $S$ has speciality $1$ and it
contains a canonical curve as a proper bisecant curve. We call $S$ a {\it canonical scroll}. 

Note that the existence of canonical scrolls is related to the existence of canonical curves having an
involution or finite morphism of degree $2$, $\gamma:C\lrw X$. The results in $\S 5$ about the
projective normality of the canonical scroll will allow to give a nice characterization of the ideal of
these canonical curves in \cite{fuentes4}.

The paper is organized in the following way:

In $\S 1$, we study the double covers of a smooth variety. We see how we can build a double cover
$\gamma:C\lrw X$ of a smooth variety $X$. We characterize when the variety $C$ is smooth. Moreover we
study the ruled variety generated by the involution on $C$. Although we will apply these results to the
case of curves we will work over smooth varieties of arbitrary dimension to give them more generality.

In $\S 2$, we use the results of the first section to prove the Segre Theorem and give the
geometrical  model of the canonical scrolls. Given a nonspecial divisor $\b$ of degree $\pi-1=deg(\b)\geq
2g-2$, we call {\it canonical geometrically ruled surface} to the ruled surface $\P(\E_{\bb})$,
$\E_{\bb}=\Te_X\oplus
\Te_X(\k-\b)$, such that the generic curve $C$ of $|2X_1|$ is smooth. When the curve $C$ is
nonhyperelliptic the image of $\P(\E_{\bb})$ by the map defined by the linear system $|X_0+\b f|$ is a
canonical scroll. From this, we conclude that any special scroll has a special directrix curve. This
last result was proved by Segre in \cite{segre1} with a condition over the degree of the scroll.
Furthermore, in this section  we see how a ruled surface is transformed by projecting from a point of a
bisecant curve.

In $\S 3$, we study  when the smooth curve $\C\in |2X_1|$ is not hyperelliptic and the complete linear
system defined by $H\sim X_0+\b f\sim X_1+\k f$ is base-point-free. This clarifies the equivalence
between both concepts: canonical scroll and canonical geometrically ruled surface. Furthermore, we
characterize the hyperelliptic double cover of smooth curves.

In $\S 4$ we prove the existence of canonical geometrically ruled surfaces over a smooth curve $X$ of
genus $g\geq 1$. 

First of all, Proposition \ref{cuandoesliso} characterizes the nonspecial divisors $\b$ such
that
$\P(\E_{\bb})$ is canonical. Such divisors satisfy the semicontinuity property that: if $P\in
X$ is a base point of the linear system $|2(\b-\k)|$, $P$ is not a base point of
$|2(\b-\k)-P|$. Then for any nonspecial divisor $\b$ such that $\deg(2\b)\geq 6(g-1)+1$, the
geometrically ruled surface $\P(\E_{\bb})$ is canonical. The proof of the existence is reduced
to the range
$\deg(\b)\leq 3(g-1)$.

In the range A: $5(g-1)\leq \deg(2\b)\leq 6(g-1)$; if $\b$ is a generic
nonspecial divisor, then $\P(\E_{\bb})$ is a canonical ruled surface, $|X_1|$ consists of a
unique curve and $2(\b-\k)$ is nonspecial. Moreover if $\b$ and $P\in X$ are generic in the
range $5(g-1)+1< \gr(2\b)\leq 6(g-1)$; the elementary transformation of $\P(\E_{\bb})$ at the
point $X_{\bb}\cap Pf$ is the general canonical ruled surface $\P(\E_{\b-P})$ in the case
$\gr(\b)-1$. 

In the range B: $4(g-1)\leq \deg(2\b)\leq 5(g-1)$; Proposition \ref{cuandoesliso} implies that
the linear system $|2(\b-\k)|$ consists of a unique divisor formed by different points, and
$\b$ cannot be a generic nonspecial divisor. Moreover, if $a:=\deg(2(\b-\k))$, let $Z_a=\{\b\in
X^b/
\P(\E_{\bb})\mbox{ is canonical} \} \subset X^b$. For $a<g-1$, if $\b\in Z_a$ is generic,
$2\b-2\k\sim P_1+\ldots+ P_a$, where all $P_i$ are different. Then, $\P(\E_{\bb})$ is the
elementary transformation of $\P(\E_{\b+P})$ but, since $2\b+2P-2\k\sim P_1+\ldots +P_a+2P$,
$\P(\E_{\b+P})$ is not canonical. That is; in range B, the generic component corresponding to
canonical geometrically ruled surfaces in case $\deg(\b)$ does not dominate the generic
component in the case $\deg(\b)-1$.

This result makes interesting to classify canonical geometrically ruled surfaces in range B. By
the existence theorem this classification is equivalent to the existence of nonspecial curves
that are not projectively generic, see Remark \ref{particular}. We hope to study their geometrical
characterization and classification  in a future paper.

Finally, in $\S 5$ we study the analogous of Noether Theorem about the projective normality of
the canonical scrolls. Theorem \ref{principal} says that the canonical scroll $R$ is
projectively normal iff the directrix curves $\ov{X_0}$ and $\ov{X_1}$ are projectively normal.
Moreover, the speciality respect to hypersurfaces of degree $m$ of $R$ is the sum of the corresponding
specialities respect to hypersurfaces of degree $m$ of the two minimal directrix curves $\ov{X_0}$ and
$\ov{X_1}$ (except when $m\geq 3$, and $g=2,deg(\b)=3$ or $g=3,deg(\b)=4$).

Therefore, in the hyperelliptic case, the canonical scroll is not projectively normal because Noether
Theorem for canonical curves; but in the nonhyperelliptic case  $R$ is projectively normal
iff the nonspecial directrix $X_{\bb}$ is projectively normal. In particular, if $\deg(\b)\geq 2g+1$,
Castelnuovo-Mumford Lemma allows to assert that the canonical scroll is projectively normal. In cases
$\deg(\b)=2g+1-k$, $k=1,2,3$, the results due to Green-Lazarsfeld \cite{green} conclude Theorem
\ref{normalidadcanonica} about the projective normality of the canonical scroll.

We thanks Lawrence Ein by his interest on this work during his visit to our Department on
November, 2001.

\bigskip

\section{Double covers of a smooth variety.}

\begin{defin}
Let $C,X$ be smooth varieties. Let $\gamma:C\lrw X$ be a surjective finite $2:1$ map. We say that
$\gamma$ is a double cover of $X$.
\end{defin}

Let $\E$ be a locally free sheaf of rank $2$ over a smooth variety $X$. Let $\pi:V=\P(\E)\lrw X$ the
corresponding ruled variety and $H$ the hyperplane section such that $\Te(1)\cong \Te_V(H)$. The fibers
of $\pi$ are isomorphic to $\P^1$. We say that a divisor $C$ on $V$ is a {\it bisecant divisor} if
$C.\pi^{-1}(x)=2$.

\begin{prop}\label{tsegre}
Let $j:C\lrw V$ be a nonsingular bisecant irreducible divisor in $V$. Let $\gamma:C\lrw X$ be the map
$\gamma=p\circ j$. Then $\E\cong \gamma_*\Te_C(H)$, $j^*\Te_V(H)\cong \Te_C(H)$ and $H^i(\E)\cong
H^i(\Te_C(H))$, for
$i\geq 0$.
\end{prop}
{\bf Proof:} Consider the exact sequence:
$$
0\lrw \Te_V(H-C)\lrw \Te_V(H)\lrw j_*\Te_C(H)\lrw 0
$$
Taking $\pi_*$, we have:
$$
0\lrw \pi_*\Te_V(H-C)\lrw \pi_*\Te_V(H)\lrw \pi_*j_*\Te_C(H)\lrw {\cal R}^1 \pi_*\Te_V(H-C)
$$

Let $x\in X$ and let $V_x\cong \P^1$ be the fiber of $\pi$ over $x$. Then:
$$
h^i(V_x,\Te_V(H-C)|_{V_x})=h^i(\Te_{P^1}(-1))=0\mbox{ when $i\geq 0$}
$$
By the Theorem of Grauert (\cite{hartshorne},III,12.9), ${\cal R}^i \pi_*\Te_V(H-C)=0$ when $i\geq 0$
and then:
$$
\E\cong \pi_*\Te_V(H)\cong \pi_*j_*\Te_C(H)=\gamma_*\Te_C(H)
$$
Furthermore, the fiber of $\gamma$ has dimension $0$ so $h^i(C_x,\Te_C(H))=0$ when $i>0$ and ${\cal
R}^i\gamma_*\Te_C(H)=0$. By (\cite{hartshorne},III,Ex 8.1), it follows that $H^i(\Te_C(H))\cong
H^i(\gamma_*\Te_C(H))$ for all $i\geq 0$.

\begin{prop}\label{morfismofinito}

Let $\gamma:\C\lrw X$ be a double cover of $X$. Let $\E_0=\gamma_*\Te_C$, then:

\begin{enumerate}

\item $\E_0$ is a locally free sheaf on $X$ of rank $2$. From this, $\pi:\P(\E_0)\lrw X$ is a
geometrically ruled variety.

\item $\E_0$ is decomposable. Moreover, $\E_0\cong \Te_X\oplus \Te_X(\eh)/ E\in Div(X)$. 

\item There is a closed immersion $j:\C\lrw \P(\E_0)$ verifying $\pi \circ j=\gamma$.

\end{enumerate}

\end{prop}

{\bf Proof:} \begin{enumerate}

\item Since $\gamma$ is a $2:1$ finite morphism:
$$
dim_{k(x)}H^0(\C_x,{\Te_{\C}}_{x})=2
$$ is constant on $X$. By Theorem of Grauert (\cite{hartshorne},III,12.9), $\E_0=\gamma_*\Te_X$ is
locally free on $X$ of rank $2$ and $\pi:\P(\E_0)\lrw X$ is a geometrically ruled variety.

\item There is a natural map:
$$
\gamma^{\#}:\Te_X\lrw \gamma_*\Te_{\C}
$$
Because $\gamma$ is dominant, $\gamma^{\#}$ is injective. Let ${\L}=Coker(\gamma^{\#})$. At any
point $x\in X$, we know that $dim(\Te_X\otimes k(x))=1$ and $dim(\gamma_*\Te_{\C}\otimes k(x))=2$. Since
$\gamma^{\#}_x$ is injective, $dim({\L}\otimes k(x))=1$ so $\L$ is an invertible sheaf on $X$.

Thus, we have the following exact sequence:
$$
0\lrw \Te_X\lrw \gamma_*\Te_{\C}\lrw \L\lrw 0
$$
We build a retraction $\phi:\gamma_*\Te_{\C}\lrw \Te_X$ of $\gamma^{\#}$. Let
$U$ be an open set of $X$ and $V=\gamma^{-1}(U)$. Given $\sigma\in \Te_{\C}(V)$ we define:
$$
\phi(U)(\sigma)(x)=\frac{1}{2}(\sigma(y_1)+\sigma(y_2))\mbox{ with }\gamma^{-1}(x)=\{y_1,y_2\}
$$
We see that $\phi\circ \gamma^{\#}=Id$, so $\phi$ is a retraction, the sequence is split exact and
$\E_0\cong \Te_X\oplus \L$ is decomposable.

\item We have the natural surjective morphism $\gamma^*\gamma_*\Te_C\lrw \Te_C$. By Propo\-sition $7.2$,
\cite{hartshorne}, this is equivalent to have a map $j:\C\lrw \P(\E_0)$ verifying $\pi
\circ j=\gamma$. We have to prove that it is a closed immersion. It is sufficient to check it in each
fibre. But, the fibres are the two points over each point of $X$ and
$\Te_C$ is very ample on these points. \qed

\end{enumerate}

\begin{rem}

{\em The decomposable ruled variety has two canonical sections corresponding to the surjections:
$$
\begin{array}{l}
{\Te_X\oplus \Te_X(E)\lrw \Te_X(E)\lrw 0}\\
{\Te_X\oplus \Te_X(E)\lrw \Te_X\lrw 0}\\
\end{array}
$$
We will denote them by  $X_0$ and $X_1$ respectively. Moreover, $X_1\sim X_0-\pi^*E$.}\qed
\end{rem}

\begin{teo}\label{canonica}
Let $\gamma:\C\lrw X$ be a double cover. Let $\E=\gamma_*\omega_C$.
Then:
\begin{enumerate}

\item $\E\cong \omega_X\oplus \Te_X(\bh)\cong \gamma_*\Te_C\otimes \Te_X(\bh)$, with $\bh\sim
\K_X-\eh$. From this, $\P(\E)\cong \P(\E_0)$. 

\item $C\sim 2X_1$ in $\P(\E)$.

\item $\K_C\sim \gamma^*\bh$ and $\eh$ is a divisor verifying $-2\eh\sim \B$, where $\B$ is the branch
divisor of
$\gamma$.

\end{enumerate}
\end{teo}
{\bf Proof:} \begin{enumerate}

\item By the Proposition \ref{morfismofinito} we know that $\gamma_*\Te_C$ is the locally free sheaf
$\E_0\cong\Te_X\oplus \Te_X(\eh)$.

By duality:
$$
\begin{array}{rl}
{\gamma_*\omega_C}&{\cong Hom(\gamma_*\Te_C,\omega_X)\cong (\gamma_*\Te_C)^*\otimes \omega_X\cong}\\
{}&{\cong (\Te_X\oplus \Te_X(\eh))^*\otimes \omega_X\cong \omega_X\oplus \Te_X(K_X-E)}\\
\end{array}
$$ 

\item $C$ is a bisecant variety on $V=\P(\E)$, that is, meets each fibre at two points. Then
$\Te_V(C)\cong \Te_C(2X_0)\oplus \pi^*(\L)$, where $\L$ is an invertible sheaf on $X$.

We use the adjunction formula:
$$
\omega_C\cong \omega_V\otimes \Te_V(C)\otimes \Te_C
$$
We know that $\omega_V\cong \Te_V(-2X_0)\otimes \pi^*(\omega_X\otimes \Te_X(\eh))$. From this,
$\omega_C\cong \Te_C\otimes \pi^*(\omega_X\otimes \Te_X(\eh)\otimes \L)$. Applying $\gamma_*$ we have:
$$
\begin{array}{rl}
{\gamma_*\omega_C}&{\cong \gamma_*(\Te_C\otimes \pi^*(\omega_X\otimes \Te_X(\eh)\otimes \L))\cong }\\
{}&{\cong \gamma_*\Te_C\otimes \omega_X\otimes \Te_X(\eh)\otimes \L}\\
\end{array}
$$
We saw that $\gamma_*\omega_C\cong \gamma_*\Te_C\otimes \Te_X(K_X-\eh)$, so $\L\cong \Te_X(-2\eh)$
and the conclusion follows.

\item We have seen that $$\omega_C\cong \Te_C\otimes \pi^*(\omega_X\otimes \Te_X(-\eh))$$ then,
 we have that $\K_C\sim \gamma^*\K_X-\gamma^*\eh$. On the other hand we know that $\K_C\cong
\gamma^*\K_X+\Ram$, where $\Ram$ is the ramification divisor of $\gamma$. From this, applying $\gamma_*$ we obtain that
$\B\sim \gamma_*\Ram\sim -\gamma_*\gamma^*\eh\sim -2\eh$, where $\B$ is the branch divisor of $\gamma$.

 \qed

\end{enumerate}

Let $\pi: V=\P(\E)\lrw X$ a decomposable ruled variety over a smooth variety $X$, such that $\E\cong
\Te_X\oplus \Te_X(\eh)$, with $\eh\not\sim 0$. Let us fix a unisecant irreducible divisor $X_1\in
|X_0-\pi^*\eh|$. We have a unique nontrivial involution $\vhi:V\lrw V$ fixing the divisors $X_0$ and
$X_1$. The unique base points of the involution are the points of $X_0$ and $X_1$. Moreover the
generators are invariant by
$\vhi$.

\begin{lemma}\label{lema1}
The unique unisecant irreducible divisors of $V$ that are invariant by $\vhi$ are $X_0$ and $X_1$.
\end{lemma}
{\bf Proof:} Let $D$ be an irreducible unisecant divisor such that it is invariant by $\vhi$. Since the
generators are invariant by $\vhi$ and $D$ meets each generator at a unique point, $D$ must
be fixed. But the unique base points of the involution are in $X_0$ and $X_1$. \qed

The involution $\vhi$ induces an involution on the linear system $|D|$ of $V$. We will denote
it by $\vhi_{|D|}:|D|\lrw |D|$.

\begin{lemma}\label{lemma2}

Let $D\sim X_0+ \pi^* \ah$ an unisecant divisor on $V$. Then $|D|$ has exactly two spaces of
base points by the involution $\vhi_{|D|}$:
$$
\begin{array}{l}
{W_0=\{ X_0+\pi^* \ch/ \ch \sim \ah \}}\\
{W_1=\{ X_1+\pi^* \ch/ \ch \sim \ah+\eh \}}\\
\end{array}
$$
\end{lemma}
{\bf Proof:} The spaces $W_0$ and $W_1$ are spaces of base points, because their divisors are
composed by invariant varieties ($X_0$, $X_1$ and generators).

Since $\dim(|D|)=h^0(\Te_X(\ah))+h^0(\Te_X(\ah+\eh))$, they are of complementary dimension. Thus
either $\vhi_{|D|}$ is the identity or they are exactly the spaces of base points.

But by Lemma \ref{lema1} the unique invariant unisecant irreducible divisors are $X_0$ or
$X_1$; so, if $\vhi_{|D|}=\id$, $W_0=|D|$ or $W_1=|D|$ and the conclusion follows. \qed

\begin{prop}\label{prop1}
The linear system $|2X_1|$ has exactly two spaces of fixed \- divisors by the involution $\vhi$:
$$
\begin{array}{l}
{F_0=\langle F_0',2X_1 \rangle,\mbox{ where } F_0'=\{2X_0+\pi^*\B/ \B\sim -2\eh \} }\\
{F_1=\{X_0+X_1+\pi^*\ch/ \ch\sim -\eh \} }\\
\end{array}
$$
\end{prop}
{\bf Proof:} The spaces $F_0'$ and $F_1$ are spaces of fixed divisors. Let $D$ be a divisor of $F_0$. $D$
is in the pencil obtained joining $F_0'$ with $2X_1$; that is, $D\in \langle 2X_0+ \pi^*\B,2X_1
\rangle=L$, with $\B\sim -2\eh$.

The divisors of $L$ are bisecant divisors. The base points of these family of divisors are the
points of $\B\subset X_1$. Thus, the pencil $L$ defines an involution in the
generic generator $P f$ of $S$ which relates points in the same generator which lie in the same
bisecant divisor of $L$. 

This involution has two fixed points $X_0\cap P f$ and $X_1\cap P f$. If we consider the trace of
$L\subset H^0(\Te_S(2X_1))$ on the generic generator $H^0(\Te_{P^1}(2))$, the image is the pencil
generated by the polynomials $\{x_0^2,x_1^2\}$, corresponding to the two fixed points. Any other
polynomial $\lambda x_0^2+\mu x_1^2$ of the pencil corresponds to two different points. From this, the
involution defined by the pencil $L$ is nontrivial.

On the other hand, we have the restriction of the involution $\vhi$ to the generator $P f$. This
involution have the same two fixed points. Because $\vhi$ is not the identity, both
involutions coincide and all bisecant divisors of $L$ are invariant by $\vhi$. 

We know that $\dim(|2X_1|)=h^0(\Te_X(-\eh))+h^0(\Te_X(-2\eh))$, $\dim(F_0)=h^0(\Te_X(-2\eh))$
and $\dim(F_1)=h^0(\Te_X(-\eh))-1$. The involution $\vhi_{|2X_1|}$ either has two fixed spaces
of complementary dimension or is the identity. The spaces $F_0$ and $F_1$ are invariant and
they have complementary dimension. 

If $\vhi_{|2X_1|}$ is the identity, the restriction to $|2X_1-X_1|$ and $|2X_0-X_0|$ is the
identity too. By Lemma \ref{lemma2}, this happens when $h^0(\Te_X(-\eh))=0$ and in this
case $F_1=\emptyset$. \qed

\begin{lemma}
Let $\B$ be a divisor of $X_1$ such that $\B\sim -2\eh$. Let $L$ be the pencil $\langle
2X_0+ \pi^*\B,2X_1 \rangle$. If the generic divisor of this pencil is smooth then is
irreducible.
\end{lemma}
{\bf Proof:} Note that the base points of the pencil are exactly in $\B\subset X_1$.

Let $D\in L$ such that $D\neq 2X_0+\pi^*\B$ and $D\neq 2X_1$.  

Suppose that $D=Y+\pi^*\ah$ such that $\ah$ is an effective divisor on $X$. Then the pencil
would have base points on $\ah\subset X_0$ and this is false.

From this, if $D$ is reducible it must be $D=D_1+D_2$ with $D_i\sim X_0+\pi^*\ah_i$, $\ah_1+\ah_2\sim
-2\eh$ two unisecant irreducible divisors. Then $D_1.D_2\sim -\eh$ is an effective divisor on $X$. Since
$\eh\not\sim 0$, $D_1$ meets $D_2$ and $D$ is not smooth. \qed

\begin{teo}
Let $\B$ be a divisor of $X_1$ such that $\B\sim -2\eh$. Then, there is a pencil
of bisecant divisors in $|2X_1|$ which are invariant by the involution and meet $X_1$ at $\B$. Moreover
the generic divisor of the pencil is irreducible and smooth if and only if $\B$ is smooth.
\end{teo}
{\bf Proof:} By the Proposition \ref{prop1}, we only have to prove that the generic curve of the pencil
$L=\langle 2X_0+ \pi^*\B,2X_1 \rangle$ is irreducible and smooth.

By Bertini's Theorem the generic element of the pencil is smooth away from the base locus. The base
locus of the pencil are the points of the divisor $\B$ on $X_1$. Let $D$ be the generic divisor of $L$.
We saw that if $D$ is smooth then it is irreducible.

If $D$ has a singular point, then $D.X_1=\B$ must have a singular point.

On the contrary, suppose that $D$ is smooth, and let $x\in \B$ be a base point of the pencil. If we 
consider the trace of the pencil over the generator $\pi^{-1}(x)$, we see that any divisor of the pencil
meets $\pi^{-1}(x)$ at $x$ with multiplicity $2$. Since $D$ is smooth, $D$ is tangent to the generator
$\pi^{-1}(x)$ and it meets $X_1$ transversally. From this, $x$ is a smooth point of $R$. \qed

\begin{rem}
{\em We have excluded the case $\eh\sim 0$. Consider the ruled variety $\P(\Te_X\oplus \Te_X)$. The
linear system $|X_1|=|X_0|$ is a pencil of irreducible unisecant varieties isomorphic to $X$. Then, all
the divisors of the linear system $|2X_1|$ are reducible pairs $D_1+D_2$ with $D_i\sim X_0$. In fact this
case corresponds to the trivial cover of $X$ by two copies of itself: $\gamma:X\cup X\lrw X$.}\qed
\end{rem}

\bigskip

\section{Canonical scrolls and canonical geometrically ruled surfaces.}\label{uno}

Let $\pi:\P(\E_0)\lrw X$ be a geometrically ruled surface over a smooth curve $X$ of genus $g>0$.
Let $H\sim X_0+\b f$ be a unisecant divisor such that the complete linear system $|H|$ provides
a birational morphism $\vhi_H:\P(\E_0)\lrw S\subset \P^n$, $\dim(|H|)=n$. The study of the
scroll
$S$ is equivalent to study the polarized geometrically ruled surface
$(\p(\E_0),\Te_{P(\E_0)}(H))$ and equivalent to the study of the locally free sheaf 
 $\E\cong \pi_*\Te_{P(\E_0)}(H)$ over $X$.

Applying the Proposition \ref{tsegre} to the case of curves we recover the following result due to
Corrado Segre and mentioned in the introduction.

\begin{teo}\label{teoremasegre}
Let $(\p(\E_0),\Te_{P(\E_0)}(H))$ be the polarized geometrically ruled surface
$(\p(\E_0),\Te_{P(\E_0)}(H))$. If $C\sim 2X_0+\b f$ is a bisecant curve, then:
$$
H^0(\Te_{P(\E_0)}(H))\simeq H^0(\Te_{C}(H));\quad H^1(\Te_{P(\E_0)}(H))\simeq
H^1(\Te_{C}(H))
$$
\end{teo} \qed

By using this result we can give the following definition.

\begin{defin}\label{scrollcanoncia}
Let $X$ be a smooth curve of genus $g>0$ and let $\C$ be a
nonhyperelliptic smooth curve of genus $\pi\geq 4$. Let $\C_K\subset P^{\pi-1}$ the canonical model of
$C$. Suppose that $C$ has an  involution over $X$. Then we call canonical scroll to the variety obtained
joining the points of $\C_K$ related by the involution.
\end{defin}

In a similar way, from the results obtained in the section above we can give the following definition.

\begin{propdef}\label{cuandoesliso}
Given a smooth curve $X$ of genus $g$ and a geometrically ruled surface $\pi:\P(\E)\lrw X$ the following
conditions are equivalent:

\begin{enumerate}

\item There is a smooth irreducible bisecant curve on $\P(\E)$, $j:\C\lrw \P(\E)$, such that
$\P(\E)\cong \P(\gamma_*\Te_C)$, where $\gamma=\pi \circ j$.

\item There is a smooth irreducible curve $\C$ and a double cover $\gamma:\C\lrw X$ such that
$\P(\E)\cong \P(\gamma_*\Te_C)$.

\item $\P(\E)\cong \P(\Te_X\oplus \Te_X(\k-\b))$ where $\b$ is a nonspecial divisor of degree
$deg(\b)\geq 2g-2$ and the generic element $\C\in |2X_1|$ is smooth.

\item $\P(\E)\cong \P(\Te_X\oplus \Te_X(\k-\b))$ where $\b$ is a nonspecial divisor of degree
$deg(\b)\geq 2g-2$ and $2(\b-\k)$ is a smooth divisor.
\end{enumerate}

A geometrically ruled surfaced verifying any of these conditions is called a canonical geometrically
ruled surface. We will denote it by $S_{\bb}=\P(\E_{\bb})$ where $\E_{\bb}\cong \Te_X\oplus \Te_X(\k-\b)$
and $\b$ is a nonspecial divisor of degree $deg(\b)\geq 2g-2$ verifying that $2(\b-\k)$ is a smooth
divisor.

\end{propdef}

Thus, a canonical scroll $R$ is the image by the complete linear system $|H|=|X_0+\b f|$ of a
canonical geometrically ruled surface $\P(\E_{\b})$, where the bisecant curve $\C\in |2X_1|$ is not
hyperelliptic. In the next section, we will see that in this case the linear system $|H|$ is
base-point-free and it defines a birational map. From this we have the following geometrical
description of a canonical scroll:

\begin{teo}\label{descripcionrcanonica}
Let $R\subset \P^{\pi-1}$ be a canonical scroll of genus $g>0$. $R$ is generated
by a correspondence between a  canonical curve of genus $g$ and
a nonspecial curve of genus $g$. They are linearly normal in  disjoint spaces that generate
$\P^{\pi-1}$.  Let $\Te_X(K)$ and $\Te_X(\b)$ the invertible sheaves of these curves. Then $deg(\b)\geq
2g-2$.

Consider the normalized geometrically ruled surface
$\P(\Te_X\oplus \Te_X(\e))$ where
$\e\sim \K-\b$. Let $X_0$ and $X_1\sim X_0-\e f$ be the minimal sections.
$R$ is the image of:
$$
\phi_H:\P(\Te_X\oplus \Te_X(\e))\lrw R\subset \P^{\pi-1}
$$
where $H\sim X_0+\b f\sim X_1+\k f$. Moreover, the restriction maps to the sections $X_0$ and
$X_1$ are the morphisms defined by $\Te_X(\K)$ and $\Te_X(\b)$ respectively.

The support of the singular locus is at most $X_0\cup
X_1$. Given the involution $\gamma:\C_{\K}\lrw X$, the branch divisor $\B$
satisfies: $\B\sim 2(\b-\k)$ and $\C\sim 2X_1$.
\end{teo}

C.Segre gives in \cite{segre1} a condition over the degree of an special scroll to have a special
directrix curve. We can see know that the condition over the degree is not necessary.

\begin{teo}
A linearly normal special scroll $R\subset \P^N$ is the projection of a canonical scroll.
\end{teo}
{\bf Proof:} 
Suppose that $R$  is defined by the ruled surface $S$ and the unisecant linear system $|H|$ on $S$.
We can take a smooth bisecant curve $C$ on $S$, verifying that the linear system $H\cap C$ on $C$ is
very ample (for example we can take $C\in |2H|$). Therefore, we have a double cover  $\gamma:C\lrw X$.
Let $\ov{C}$ the image of $C$ by the map defined by the linear system $|H|$ on $C$.

By Theorem \ref{tsegre} we know that $\ov{C}$ and $R$ have the same speciality. The curve $C$ is not
hyperelliptic, because an special divisor over an hyperelliptic curve is never very ample. Thus,
$C$ is the projection of a canonical curve
$C_K$ of genus $g$
$g$, and then the scroll $R$ is the projection of the canonical scroll $R_{\b}$ defined by
$\gamma$ over $C_K$. \qed

\begin{cor}
A linearly normal special scroll $R\subset \P^N$ always has an special directrix curve.
\end{cor}
{\bf Proof:} We saw that a canonical scroll
$R_{\b}$ has a canonical special directrix curve. This curve goes to an special curve for any projection
of $R_{\b}$. Now, it is sufficient to apply the above theorem. \qed

Finally, we see how a ruled surface is transformed by projecting from a point of a bisecant curve.

\begin{prop}\label{proycanonica}
Let $\pi:S\lrw X$ a ruled surface. Let $C\subset S$ be a bisecant smooth curve
 and $\gamma:C\lrw X$ the corresponding double cover. Let $\aa$ be a divisor on $C$ and $H$ an
unisecant divisor on
$S$ such that $\gamma_*\Te_C(\aa)\cong \pi_*\Te_S(H)$, or equivalently, $\Te_C(H)\cong \Te_C(\aa)$.
Then, we have that:
\begin{enumerate}

\item If $\b$ is a divisor on $X$:
$$
\gamma_*\Te_C(\aa+\gamma^*\b)\cong \pi_*\Te_S(H+\b f)
$$

\item If $x$ is a point of $C$, such that $\gamma(x)=P$:
$$
{\gamma_*\Te_C(\aa-x)\cong \pi'_*\Te_{S'}(\nu^*(H)-Pf)}\\
$$
where $\nu:S'\lrw S$ is the elementary transform of $S$ in $x$.

\item If $x$ is a point of $C$, such that $\gamma(x)=P$:
$$
{\gamma_*\Te_C(\aa+x)\cong \pi'_*\Te_{S'}(\nu^*(H))}\\
$$
where $\nu:S'\lrw S$ is the elementary transform of $S$ in $\gamma^*(P)-x$.

\end{enumerate}
\end{prop}
{\bf Proof:} \begin{enumerate}

\item Let $\b$ be a divisor on $X$. We use the projection formula:
$$
\begin{array}{rl}
{\gamma_*\Te_C(\aa+\gamma^*\b)}&{\cong \gamma_*(\Te_C(\a)\otimes \gamma^*\Te_X(\b))\cong
\gamma_*\Te_C(\aa)\otimes \Te_X(\b)\cong}\\
{}&{\cong \pi_*\Te_S(H)\otimes \Te_X(\b)\cong \pi_*\Te_S(H+\b f)}\\
\end{array}
$$

\item Let $x$ be a point of $C$ such that $\gamma(x)=P$. Consider the exact sequence:
$$
0\lrw \Te_C(\aa-x)\lrw \Te_C(\aa)\lrw \Te_x\lrw 0
$$
Applying $\gamma_*$ and because ${\cal R}^1\gamma_*\Te_C(\aa-x)=0$, we obtain:
$$
0\lrw \gamma_*\Te_C(\aa-x)\lrw \gamma_*\Te_C(\aa)\lrw \Te_P\lrw 0
$$
In this way, we see that $\P(\gamma_*\Te_C(\aa-x))=S'$ is the elementary transform of $S$ in the point
$x$.Moreover, if we denote by
$\nu$ the elementary transformation,
$$
\gamma_*\Te_C(\aa-x)\cong \pi'_*\Te_{S'}(\nu^*(H)-Pf)
$$

\item Let $x$ be a point of $C$ such that $\gamma(x)=P$. Let $y=\gamma^*(P)-x$. We have that:
$$
\gamma_*\Te_C(\aa+x)\cong \gamma_*\Te_C(\aa+\gamma^*(P)-y)
$$
But, applying the case $1$ of this Proposition:
$$
\gamma_*\Te_C(\aa+\gamma^*(P))\cong \pi_*\Te_S(H+Pf)
$$
and now, by using the above case, if $\nu:S'\lrw S$ the elementary transformation of
$S$ at the point
$y$, we see that:
$$
\gamma_*\Te_C(\aa+\gamma^*(P)-y)\cong \pi'_*\Te_{S'}(\nu^*(H+Pf)-Pf)\cong \pi'_*\Te_{S'}(\nu^*(H))
$$
\qed
\end{enumerate}

\bigskip

\section{The elliptic and the hyperelliptic cases.}\label{unoymedio}

\begin{prop}\label{h1}
Let $\C$ be an hyperelliptic curve of genus $\pi\geq 2$ and $\gamma:\C\lrw X$ an involution of genus
$g\geq 1$. Then $X$ is elliptic or hyperelliptic. Moreover, the $g_2^1$ of $X$ makes the following
diagram commutative:
$$
\setlength{\unitlength}{5mm}
\begin{picture}(14,3.6)

\put(2.6,3){\makebox(0,0){$\C$}}
\put(2.6,0){\makebox(0,0){$X$}}
\put(10,3){\makebox(0,0){$\P^1$}}
\put(10,0){\makebox(0,0){$\P^1$}}
\put(10,2.3){\vector(0,-1){2}}
\put(2.6,2.3){\vector(0,-1){2}}
\put(3,3){\vector(1,0){6}}
\put(3,0){\vector(1,0){6}}
\put(6.3,3.6){\makebox(0,0){$\pi_2^1$}}
\put(6.3,0.6){\makebox(0,0){$g_2^1$}}
\put(1.9,1.5){\makebox(0,0){$\gamma$}}
\put(10.9,1.5){\makebox(0,0){$\ov{\vhi}$}}

\end{picture}
$$
where $\ov{\vhi}$ is a $2:1$ morphism and $\pi_2^1$ is the pencil of degree $2$ of $\C$. 
\end{prop}
{\bf Proof:} Let $\mu: C\lrw C$ the automorphism defined by the involution $\gamma$. Because $\C$ is
hyperelliptic, it has a unique $\pi_2^1$. Therefore, there exist an isomorphism $\vhi:\P^1\lrw \P^1$
verifying $\pi_2^1\circ \mu=\vhi\circ \pi_2^1$.

$\vhi$ induces a $2:1$ morphism $\ov{\vhi}:\P^1\lrw \P^1$ parameterizing the points related by $\vhi$,
that is, $\ov{\vhi}(x)=\ov{\vhi}(y)\Leftrightarrow \vhi(x)=y$. Then, we have:
$$
\setlength{\unitlength}{5mm}
\begin{picture}(14,3.6)

\put(2.6,3){\makebox(0,0){$\C$}}
\put(2.6,0){\makebox(0,0){$X$}}
\put(10,3){\makebox(0,0){$\P^1$}}
\put(10,0){\makebox(0,0){$\P^1$}}
\put(10,2.3){\vector(0,-1){2}}
\put(2.6,2.3){\vector(0,-1){2}}
\put(3,3){\vector(1,0){6}}
\put(3,0){\vector(1,0){6}}
\put(6.3,3.6){\makebox(0,0){$\pi_2^1$}}
\put(6.3,0.6){\makebox(0,0){$g_2^1$}}
\put(1.9,1.5){\makebox(0,0){$\gamma$}}
\put(10.9,1.5){\makebox(0,0){$\ov{\vhi}$}}

\end{picture}
$$
where $g_2^1(P)=\ov{\vhi}(\pi_2^1(\gamma^{-1}(P)))$ and it is well defined, because:
$$
\gamma^{-1}(P)=\{Q,\mu(Q)\} \mbox{ and } \vhi(\pi_2^1(Q))=\pi_2^1(\mu(Q)) \Rightarrow
\ov{\vhi}(\pi_2^1(\mu(Q)))=\ov{\vhi}(\pi_2^1(Q)). 
$$
Thus $X$ has a $g_2^1$ and is elliptic or hyperelliptic. \qed 

\begin{cor}\label{h2}
If $X$ is neither elliptic nor hyperelliptic the linear system $H\sim X_0+\b f$ is base-point-free and
it defines a birational map.
\end{cor}
{\bf Proof:} Because $|\K|$ is very ample, it is sufficient to see that the linear
system $|\b|$ is base-point-free:
\begin{enumerate}

\item If $deg(\b)\geq 2g$, the divisor $\b$ is always base-point-free.

\item If $deg(\b)=2g-1$ and $\b$ has a base point, then $\b\sim \k+P$. Since $\P(\E_{\b})$ is a
canonical ruled surface, $2(\b-\k)\sim 2P$ is linearly equivalent to two different points. Then $X$ has
a $g_2^1$, but this is false by hypothesis. 

\item If $deg(\b)=2g-2$ and $\b$ has a base point, then $\b-P\sim \k-Q$. But $2(\b-\k)\sim 2P-2Q$ must
be  linearly equivalent to divisor $0$. Then $X$ has a $g_2^1$, and this is false by hypothesis. \qed

\end{enumerate}

\begin{prop}\label{genero3}
A non-hyperelliptic curve $\C$ of genus $\pi=3$ does not have involutions of genus $2$.
\end{prop}
{\bf Proof:} Let $\gamma:\C\lrw X$ an involution of genus $2$ of $\C$. We know that $\C\in |2X_1|\subset
S$ with $S=\P(\E_{\bb})$ and $\b$ is a nonspecial divisor of degree $2$. Thus, if $\b=P+Q$ then
$h^0(\Te_X(\b))=1=h^0(\Te_X(\b-P))$. From this:
$$
h^0(\Te_C(\K_C-\gamma^*(P)))=h^0(\Te_S(H-P f))=h^0(\Te_S(H))-1=h^0(\Te_C(\K_C))-1
$$
so the canonical divisor of $\C$ is not very ample and $\C$ is hyperelliptic. \qed

\begin{prop}
Let $X$ be an elliptic curve and $S=\P(\E_{\bb})$ a canonical ruled surface with $deg(\b)=2$. Let $\C$
be a smooth irreducible curve of genus $3$. We have a double cover
$\gamma:\C\lrw X$ with branch divisor $\R=a_1+a_2+a_3+a_4$. Then $\C$ is hyperelliptic if and only if
$a_1+a_2\sim a_3+a_4\sim \b$.
\end{prop}
{\bf Proof:} Let $\Ram=A_1+A_2+A_3+A_4$ be the ramification divisor such that $\gamma(A_i)=a_i$. We have
that $\K_C\sim \gamma^*\K_X+\Ram$.

 Suppose that $\C$ is hyperelliptic. Then we have that $A_1+A_2\sim
A_3+A_4\sim
\pi_2^1$. Applying
$\gamma$ we see that $a_1+a_2\sim a_3+a_4\sim \gamma_*{\pi_2^1}$. Moreover:
$$
\begin{array}{l}
{h^0(\Te_S(H-\gamma_*{\pi_2^1}f))=h^0(\Te_C(\K_C-\gamma^{*}\gamma_*{\pi_2^1}))=1}\\
{}\\
{h^0(\Te_S(H-\gamma_*{\pi_2^1}f))=h^0(\Te_X(-\gamma_*{\pi_2^1}))+h^0(\Te_X(\b-\gamma_*{\pi_2^1}))}\\
\end{array}
$$
and it follows that $\gamma_*{\pi_2^1}\sim \b$. 

On the contrary suppose that $a_1+a_2\sim a_3+a_3\sim \b$. Then $2A_1+2A_2\sim \gamma^*(a_1+a_2)\sim
\gamma^*\b\sim K_C \sim A_1+A_2+A_3+A_4$. From this $A_3+A_4\sim A_1+A_2$ and $\C$ has a $\pi_2^1$. \qed

\begin{prop}\label{h3}
Let $S=\P(\E_{\bb})$ be a canonical ruled surface and $\C\sim 2X_1$ a smooth curve. Then $\C$ is
hyperelliptic if and only if:
\begin{enumerate}
\item $X$ is elliptic and $deg(\b)=1$ or $deg(\b)=2$ and the branch divisor $\B=a_1+a_2+a_3+a_4$
verifies $a_1+a_2\sim a_3+a_4\sim \b$.

\item $X$ is hyperelliptic and $h^0(\Te_X(\b-g_2^1))=h^0(\Te_X(\b))-1$.
\end{enumerate}
Moreover, $\C$ is elliptic if and only if $X$ is elliptic and $deg(\b)=0$.
\end{prop}
{\bf Proof:} Let us suppose that $\C$ is hyperelliptic. By Proposition \ref{h1} we know that $X$ is
elliptic or hyperelliptic and we have that $\gamma^*(g_2^1)\sim 2\pi_2^1$. From this:
$$
\Te_{\C}(H-g_2^1 f)\sim \Te_{\C}(\K_{C}-2\pi_2^1)
$$
By Theorem \ref{teoremasegre}:
$$
h^0(\Te_S(H-g_2^1 f))=h^0(\Te_{\C}(\k_{C}-2\pi_2^1))=h^0(\Te_{\C}(\k_{C})-2=h^0(\Te_S(H))-2
$$
But,
$$
\begin{array}{l}
{h^0(\Te_S(H-g_2^1 f))=h^0(\Te_X(\b-g_2^1))+h^0(\Te_X(\k-g_2^1))=}\\
{=h^0(\Te_X(\b-g_2^1))+h^0(\Te_X(\k))-1}\\
{}\\
{h^0(\Te_S(H))=h^0(\Te_X(\b))+h^0(\Te_X(\k))}\\
\end{array}
$$
and we obtain $h^0(\Te_X(\b-g_2^1))=h^0(\Te_X(\b))-1$. If $X$ is elliptic this happens when
$deg(\b)=1,2$. Moreover is $deg(\b)=2$ the conditions of Proposition \ref{genero3} holds.

Conversely, if we suppose that $X$ is hyperelliptic and $h^0(\Te_X(\b-g_2^1))=h^0(\Te_X(\b))-1$.
We can check that:
$$
h^0(\Te_{\C}(\K_{\C}-\gamma^*(g_2^1)))=h^0(\Te_S(H-g_2^1f))=h^0(\Te_S(H))-2=h^0(\Te_{\C}(\K_{\C}))-2
$$
Thus, $\C$ has a $g_4^2$. If the genus of $\C$ is $\pi\geq 4$ then $\C$ is hyperelliptic. Moreover, we
have seen that a curve of genus $3$ with an involution of genus $2$ is hyperelliptic (Proposition
\ref{genero3}). 

If $X$ is elliptic and $deg(\b)=1$ then the genus $\pi$ of $\C$ is $2$ and $\C$ is hyperelliptic. If
$deg(\b)=2$ we apply Lemma \ref{genero3}.

Finally, by Hurwitz's formula, $\C$ is elliptic if and only if $X$ is elliptic and $deg(\b)=0$. \qed

\begin{teo}\label{h4}
Let $S=\P(\E_{\b})$ a canonical ruled surface and $\C\sim 2X_1$ a smooth curve. Then $\C$ is
hyperelliptic if and only if:

\begin{enumerate}

\item  $X$ is hyperelliptic and one of the following conditions holds:

\begin{enumerate}

\item $\b\sim \k+g_2^1$.

\item $\b\sim \k+P$, with $2P\sim g_2^1$.

\item $\b\sim \sum_1^{g-2} g_2^1+P+Q$, with $2P\sim 2Q\sim g_2^1$.

\end{enumerate}

\item $X$ is elliptic and $deg(\b)=1$ or $deg(\b)=2$ and the branch divisor $\B=a_1+a_2+a_3+a_4$
verifies $a_1+a_2\sim a_3+a_4\sim \b$.

\end{enumerate}

Moreover, the cases $(b)$,$(c)$ or $deg(\b)=0,1$ are the unique cases where $\b$ has base-points.
\end{teo}
{\bf Proof:} Suppose that $\C$ is hyperelliptic. We apply  Proposition \ref{h3}. Let us study the
hyperelliptic case. If
$deg(\b)\geq 2g+1$ then
$h^0(\Te_X(\b-g_2^1))=h^0(\Te_X(\b))-2$. It remains to check the following cases:

\begin{enumerate}

\item If $deg(\b)=2g$ and $h^0(\Te_X(\b-g_2^1))=h^0(\Te_X(\b))-1$ then $\b-g_2^1$ is special of degree
$2g-2$, that is, $\b\sim \k+g_2^1$.

\item If $deg(\b)=2g-1$ and $h^0(\Te_X(\b-g_2^1))=h^0(\Te_X(\b))-1$ then $\b-g_2^1$ is special of degree
$2g-3$, that is, $\b-g_2^1\sim \k-P' \Rightarrow \b\sim \k+g_2^1-P'\sim \k+P$ with $P+P'\sim g_2^1$.
Since $\b$ defines a canonical ruled surface, $2(\b-\k)\sim 2P'$ must be linearly equivalent to two
different points, so $2P'\sim g^1_2$.

\item If $deg(\b)=2g-2$ and $h^0(\Te_X(\b-g_2^1))=h^0(\Te_X(\b))-1$ then $\b-g_2^1$ has speciality $1$ 
and degree $2g-4$, that is, $\b-g_2^1\sim \k-P-Q$ with $P+Q\not\sim g_2^1$. Because $2(\b-\k)$ must
be effective, $2P+2Q\sim 2g_2^1\Rightarrow 2P\sim 2Q\sim g_2^1$ and from this $\b\sim
\sum_1^{g-2} g_2^1+P+Q$.

\end{enumerate}

Conversely, if $X$ is hyperelliptic and one of the three conditions holds, then 
$h^0(\Te_X(\b-g_2^1))=h^0(\Te_X(\b))-1$ and by Proposition \ref{h3} the curve $\C$ is hyperelliptic.

Finally, let us suppose $\b$ has a base point. Then $deg(\b)\leq 2g-1$ and by Corollary \ref{h2}, $X$ is
elliptic or hyperelliptic. If $X$ is elliptic, a nonspecial divisor of degree $\geq 0$ has a base point
if and only if $deg(\b)=0,1$ and $\b\not\sim 0$. Suppose that
$X$ is hyperelliptic. If
$deg(\b)=2g-1$ then $\b \sim
\k+P$ and we are in the case
$(b)$. If
$deg(\b)=2g-2$ then $\b-P\sim \k-Q\Leftrightarrow \b\sim\k+P-Q$. Since $2(b-\k)\sim 0$, $2P\sim 2Q\sim
g_2^1$ and
$\b\sim \k+P-Q\sim \sum_1^{g-2} g_2^1+P+Q$ (case $(c)$). \qed

\begin{cor}\label{h5}
Let $S=\P(\E_{\b})$ a canonical ruled surface and $H\sim X_0+\b f$. The linear system $H$ defines a
canonical scroll except when: 

\begin{enumerate}

\item  $X$ is hyperelliptic and one of the following conditions holds:

\begin{enumerate}

\item $\b\sim \k+g_2^1$.

\item $\b\sim \k+P$, with $2P\sim g_2^1$.

\item $\b\sim \sum_1^{g-2}g_2^1+P+Q$, with $2P\sim 2Q\sim g_2^1$. 

\end{enumerate}

\item $X$ is elliptic and $deg(\b)=0,1,2$. \qed

\end{enumerate}

\end{cor}

\begin{cor}
The unique possible involutions of an hyperelliptic curve of ge\-nus $\pi$ are of genus
$\frac{\pi-1}{2},\frac{\pi}{2}$, or $\frac{\pi+1}{2}$. \qed
\end{cor}

\begin{cor}
If $\P(\E_{\bb})$ is a canonical ruled surface and $H\sim X_0+\b f$ defines a canonical scroll then $|H|$
is base-point-free and it defines a birational map.
\end{cor}
{\bf Proof:} If $X$ is nonhyperelliptic it is Corollary \ref{h2}.

If $X$ elliptic or hyperelliptic and $\b$ defines a canonical scroll. We saw at Proposition
 \ref{h3} that $h^0(\Te_X(\b-g_2^1))=h^0(\Te_X(\b))-2$. $|K|$ and $|\B|$ are base-point-free.  
Moreover, $|\k-P|$ and $|\b-P|$ don't have common base points for all $P\in X$. From this $|H|$ is
base-point-free and it defines a birational map. \qed

We finish this section by studying the map defined by the linear system $|\b|$ when $X$ is
hyperelliptic.

\begin{prop}\label{descripcionb}
Let $\b$ be a nonspecial divisor of degree $b$ over an hyperelliptic curve $X$ of genus $g$, with $g\geq
2$ and $b\geq 2g-2$. Suppose that $\b$ defines a canonical scroll. Let $\phi_{\b}:X\lrw \P^{b-g}$
be the map defined by the linear system $|\b|$. We have:
\begin{enumerate}

\item $\phi_{\b}$ is birational except when $g=2$ and $b=3$ or $g=3$ and $b=4$.

\item $\phi_{\b}$ is an isomorphism ($\b$ is very ample) except when:

\begin{enumerate}

\item $b=2g$ and 

\begin{enumerate}

\item $g=2,3$ and $\b\sim \k+P_1+P_2$, or

\item $g>3$ and $\b\sim \k+P_1+P_2$, with $P_1,P_2$ ramification points of the $g_2^1$.

\end{enumerate}

\item $b=2g-1$ and 

\begin{enumerate}

\item $g=2,3$, or

\item $g>3$ and $\b\sim \sum_1^{g-2} g_2^1+P_1+P_2+P_3$, with $P_1,P_2,P_3$ ramification points of the
$g_2^1$. $\phi_{\b}$ has a triple point which is its unique singular point.

\end{enumerate}

\item $b=2g-2$ and 

\begin{enumerate}

\item $g=3$, or

\item $g>3$ and $\b\sim \sum_1^{g-3} g_2^1+P_1+P_2+P_3+P_4$, with $P_1,P_2,P_3,P_4$ ramification points
of the
$g_2^1$. $\phi_{\b}$ has a
quadruple point which is its unique singular point.

\end{enumerate}

\end{enumerate}

\end{enumerate}

\end{prop}
{\bf Proof:} By Riemann-Roch Theorem we know that $\b$ is very ample when $b\geq 2g+1$.

Moreover, $\phi_{\b}$ fails to be an isomorphism iif there are two points $P_1$, $P_2$ satisfying
$\b-P_1-P_2$ is an special divisor with speciality $1$.

If $b=2g$ then $deg(\b-P_1-P_2)=2g-2$. From this, if $\b-P_1-P_2$ is special, $\b\sim \k+P_1+P_2$.
Because
$\b$ defines a canonical ruled surface $P_1+P_2\not\sim g_2^1$ and $2\b+2\k$ is linearly equivalent
to four different points. If $g=2,3$ this always happens. If $g>3$ the unique $g^1_4$ or $g^2_4$ are
composed of
$g^1_2$. Since $P_1+P_2\not\sim g_2^1$, necessary $2P_1\sim 2P_2\sim g_2^1$.

If $b=2g-1$ then $deg(\b-P_1-P_2)=2g-3$. Now, if $\b-P_1-P_2$ is special, $\b\sim \k+P_1+P_2-P_3$.
Because
$\b$ defines a canonical ruled surface $P_1+P_2\not\sim g_2^1$ and $2\b+2\k$ is linearly equivalent
to two different points. If $g=2$, $\phi_{\b}$ is a $3:1$ morphism over a line. If $g=3$, $\phi_{\b}$
projects
$X$ into a plane curve of degree
$5$ and genus
$3$ so it has a triple singular point. If $g>3$ the unique $g^1_4$ or
$g^2_4$ are composed of $g_2^1$. Since $P_1+P_2\not\sim g_2^1$, necessary $2P_1\sim 2P_2\sim
2P_3\sim g_2^1$, so $\b\sim \k+P_1+P_2-P_3\sim \sum_1^{g-2} g_2^1+P_1+P_2+P_3$ and $\phi_{\b}$ has a
triple point which is its unique singular point.

If $b=2g-2$ then $deg(\b-P_1-P_2)=2g-4$. If $\b-P_1-P_2$ is special with speciality $1$, $\b\sim
\k+P_1+P_2-P_3-P_4$ with $P_3+P_4\not\sim g^1_2$. Because
$\b$ defines a canonical ruled surface $P_1+P_2\not\sim g_2^1$ and $2\b-2\k$ is linearly equivalent
to $0$.In this case $g\neq 2$ because a nonspecial divisor of degree $2$ over a curve of genus $X$
has base points.If $g=3$, then $\phi_{\b}$ is a $4:1$ morphism over a line. If $g>3$ the unique $g^1_4$
or $g^2_4$ are composed of $g_2^1$. Since $P_1+P_2\not\sim g_2^1$, necessary $2P_1\sim 2P_2\sim
2P_3\sim 2P_4\sim g_2^1$, so $\b\sim \k+P_1+P_2-P_3-P_4\sim \sum_1^{g-3} g_2^1+P_1+P_2+P_3+P_4$ and $\phi_{\b}$ has a
quadruple point which is its unique singular point. \qed

\bigskip

\section{The existence theorem.}\label{dos}

Given a smooth curve $X$ of genus $g$, let us study the divisors $\b$ defining a canonical ruled
surface. We see that $\P(\E_{\bb})$ is a canonical ruled surface if and only if $2(\b-\k)$ is a smooth
divisor, that is, if $P$ is a base point of $2(\b-\k)$, $P$ is not a base point of $2(\b-\k)-P$. In
particular $2(\b-\k)$ is an effective divisor.

\begin{rem}\label{particular}

{\em Note that a necessary condition for $\b$ to define a canonical ruled surface is that the divisor
$2(\b-\k)$ must be an effective divisor.

Let $\a=2\b-2\k$ and $a=\deg(\a)$. If $a\geq g$ this divisor is always an effective divisor.

If $a<g$, the generic divisor of degree $a$ is not an effective divisor.  In this case, we will
see that the generic divisor $\b$ does not define a canonical geometrically ruled surface.
Moreover, if
$b\leq 3g-4$ the points of $\a$ are contained in a hyperplane section of $X_{\bb}$. If $a<g$
this inequality holds, so we have $\b\sim \a+\c$, where $\c$ is an effective divisor. Now, we
see that:
$$
\a\sim 2\b-2\k  \Rightarrow \b \sim 2\k+\a-\b \sim 2\k-\c
$$
and
$$
2\k-2\c\sim 2\k-(2\k-\b)\sim 2(\b-\k) \mbox{ is effective by hypothesis}
$$
We conclude that the hyperplane sections of $X_{\bb}$ are the residual points of the system of
quadrics sections of the canonical curve of genus $g$ passing through a divisor $\c$, such that
$\c$ are a set of contact points of a quadric with the canonical curve.

Since $h^0(\Te_X(2\k))=3(g-1)$, if $2\deg(\c)\geq 3(g-1)$ the number of conditions imposed by
$2\c$ is greater than the dimension of the linear system $|2\k|$, so the points of $\c$ cannot
be generic. This happens when $2\deg(\b)\leq 5(g-1)$, so in this range we hope that the generic
divisor $\b$ does not define a canonical ruled surface.}\qed

\end{rem}

\begin{prop}\label{libredepuntosfijosdobles}
Let $X^a$ be the family of effective divisors of degree $a$ and $U^a=\{ \a\in X^a/|\a| \mbox{
is smooth} \}$ which is an open subset of $X^a$. Then:

\begin{enumerate}

\item If $a\geq 2g-1$ then $U^a=X^a$ and any divisor is base-point-free and nonspecial.

\item If $g+1\leq a\leq 2g-2$ then $U^a\neq \emptyset$ and the generic divisor of $U^a$
is base-point-free and nonspecial. 

\item If $0\leq a\leq g$ then $U^a\neq \emptyset$ and the generic divisor $\a$ of
$U^a$ is formed by different points and it has speciality $h^1(\a)=g-a$. \qed

\end{enumerate}

\end{prop}
{\bf Proof:} From the semicontinuity of the cohomology and the Riemann-Roch Theorem we see that:

\begin{enumerate}

\item If $\deg(\a)\geq 2g$, the divisor $\a$ is always base-point-free, so it is smooth.

If $\deg(\a)=2g-1$, the divisor $\a$ has a base point $P_0$ iff $\a\sim \k+P_0$. But, in this
case $\a-P_0\sim \k$ is base-point-free, and $\a$ is smooth.

\item If $g+1\leq \deg(\a) \leq 2g-2$, the generic effective divisor $\a$ is base-point-free
and, in particular, it is smooth and nonspecial.

\item If $0\leq \deg(\a)\leq g$, the generic effective divisor is formed by different points, so
it is smooth and it has the expected speciality $h^1(\Te_X(\a))=g-deg(\a)$.
If $deg(\a)=1$ then any effective divisor $\a$ is smooth. \qed

\end{enumerate}

Given a nonspecial divisor $\b$, we will apply the results \ref{cuandoesliso},
\ref{libredepuntosfijosdobles} to study when it defines a canonical ruled surface. We will
denote by $\Te_X(\b)$ the invertible sheaves of $\Pic^b(X)$, where $\b$ is a divisor of degree
$b$. Although we will characterize the invertible sheaf $\Te_X(\b)$, the condition is
given over the divisor $2\b$ so we will use the following lemma:

\begin{lemma}\label{pordos}
Let $X$ be a smooth curve of genus $g$ and $b\geq 0$. Consider the map:
$$
\begin{array}{rcl}
{\vhi^b:\Pic^b(X)}&{\lrw}&{\Pic^{2b}(X)}\\
{\L}&{\mapsto}&{\L^2}
\end{array}
$$
Then $\vhi^b$ is a surjection with finite fibre (nontrivial when $g>0$).
\end{lemma}
{\bf Proof:} Given two invertible sheaves $\L,\L'\in \Pic^b(X)$, $\vhi^b(\L)=\vhi^b(L')$
iff $(\L^{-1}\otimes \L')^2\cong \Te_X$, that is, the fibre of $\vhi^b$ is isomorphic to the kernel of
$\vhi^0$. But $\ker(\vhi^0)= H^1_{et}(X,{\bf Z}_2)$ which is a finite group (but
nontrivial if
$g>0$).
\qed

Let $a=2b-2(g-2)$. Consider the following incidence variety and its projections:
$$
\setlength{\unitlength}{5mm}
\begin{picture}(22,4)

\put(7,2){\makebox(0,0){$J=\{(\Te_X(\b),\a)\in \Pic^b(X)\times X^a/ 2\b-2\k\sim \a\}$}}
\put(15,1.9){\vector(4,-1){4}}
\put(15,2.1){\vector(4,1){4}}
\put(20.7,3){\makebox(0,0){$\Pic^b(X)$}}
\put(20,1){\makebox(0,0){$X^a$}}
\put(17,1){\makebox(0,0){$q$}}
\put(17,3){\makebox(0,0){$p$}}

\end{picture}
$$
and 
$$
J_0=\{(\Te_X(\b),\a)\in \Pic^b(X)\times U^a/ 2\b-2\k\sim \a\}\subset J
$$
$J$ is a closed variety because it is defined by a closed condition. $J_0$ is an open set of
$J$. The invertible sheaves of $\Pic^b(X)$ that define a canonical ruled surface are $p(J_0)$,
except when $b=2g-2$. In this case, we have to discard the canonical divisor; that is, the
invertible sheaves defining a canonical ruled surface are 
$p(J_0)-\{\Te_X(\k)\}$. Because the projection maps are open and closed maps, we note
that $p(J_0)$ is an open subset of $p(J)$ and $p(J)$ is a closed subset of $\Pic^b(X)$.

$q$ is a surjection, because given $\a\in X^a$, by Lemma \ref{pordos} we know that there is
a divisor
$\b\in X^b$ such that $2\b\sim 2\k+\a$.

Moreover, given $\a\in X^a$, the generic fibre is $$q^{-1}(\a)=\{\Te_X(\b)\in \Pic^b(X)/2\b\sim
2\k+\a\}$$By Lemma \ref{pordos} we know that there are at most a finite number of invertible
sheaves verifying this condition, so $q$ has finite fibre.

Thus, $\dim(J)=a$ and by applying Proposition \ref{libredepuntosfijosdobles}, we know that
$J_0$ is a nonempty open subset of $J$, $\dim(J_0)=\dim(J)=a$ and
$\dim(p(J_0))=\dim(p(J))=a-\dim(p^{-1}(\Te_X(\b)))$.

Let us see which is the generic fibre of the map $p$. Let $\Te_X(\b)\in p(J)$ be a generic
invertible sheaf on the image. Then $2\b-2\k$ is an effective divisor. Therefore:
$$
p^{-1}(\Te_X(\b))=\{\a\in X^a/\a\sim 2\b-2\k \}=|\a|
$$

We distinguish some cases:

\begin{enumerate}

\item If $a\geq 2g-1$, the generic divisor $\a \in X^a$ is nonspecial and $U^a=X^a$. Thus,
$\dim(p(J_0))=\dim(p(J))=a-(a-g)=g$  and $p(J_0)=\Pic^b(X)$.

\item If $g\leq a \leq 2g-2$, the generic divisor $\a$ of $X^a$ is nonspecial;
$\dim(p(J_0))=\dim(p(J))=g$, where $p(J_0)$ is an open subset of $\Pic^b(X)$.

\item If $1\leq a \leq g-1$, the generic divisor $\a$ of $X^a$ has speciality $g-a$. Thus,
$\dim(p(J_0))=\dim(p(J))=a$, where $p(J)$ is a closed subset of codimension $g-a$ in $\Pic^b(X)$ and
$p(J_0)$ is an open subset of $p(J)$. In this case the generic invertible sheaf of $\Pic^b(X)$ does not
define a canonical ruled surface. 

\item If $a=0$ then $X^a=U^a$ has a unique divisor $0$. Thus
$J=J_0\cong p(J_0)=p(J)$ is a set of dimension $0$ (codimension $g$). In this case an invertible
sheaf  $\Te_X(\b)$ of $\Pic^b(X)$ defines a canonical ruled surface iff $2\b\sim 2\k$. By
Lemma \ref{pordos} they are a finite (nontrivial) number.

So, we conclude:

\end{enumerate}

\begin{teo}\label{existenciafuerte}
Let $X$ be a smooth nonhyperelliptic curve of genus  $g\geq 1$. Let
$\Pic^b(X)$ be the family of invertible sheaves $\Te_X(\b)$ of degree $b\geq
2g-2$. Then:
\begin{enumerate}

\item If $b\geq 3(g-1)+1$, any invertible sheaf of degree $b$ defines a canonical ruled surface.

\item If $(5(g-1)+1)/2 \leq b \leq 3(g-1)$, the generic invertible sheaf of degree $b$ defines a
canonical ruled surface.

\item If $2(g-1)\leq b \leq 5(g-1)/2$, the generic invertible sheaf of degree $b$ does not
define a canonical ruled surface, but there is a family of codimension $5(g-1)+1-2b$ of
invertible sheaves of degree $b$ such that the generic element defines a canonical ruled
surface.\qed

\end{enumerate}

\end{teo}

\bigskip

\section{Projective normality of the canonical scroll.}\label{normality}

In this section we will study the projective normality of the canonical scroll. Let us first remember
some definitions and results about normality of curves and decomposable ruled surfaces:

\begin{defin}
Let $V$ be a projective variety. Let ${\cal F}_i,i=1,\ldots ,s$, coherent sheaves on $V$. We
call
$s({\cal F}_1,\ldots,{\cal F}_s)$ the cokernel of the map:
$$
H^0({\cal F}_1)\otimes \ldots \otimes H^0({\cal F}_s)\lrw H^0({\cal F}_1\otimes \ldots \otimes {\cal
F}_s)
$$
If ${\cal F}_i$ are invertible sheaves $\Te_V(D_i)$ where $D_i$ are divisors on $V$, we will write
$s(D_1,\ldots,D_s)$.
\end{defin}

\begin{lemma}\label{lematecnico}
If $s({\cal F}_1,{\cal F}_2)=0$, then:
$$
s({\cal F}_1,{\cal F}_2,{\cal F}_3,\ldots,{\cal F}_s)=s({\cal F}_1\otimes {\cal
F}_2,{\cal F}_3,\ldots,{\cal F}_s)
$$
\end{lemma}
{\bf Proof:}  It sufficient note that $s({\cal F}_1,{\cal F}_2,{\cal F}_3,\ldots,{\cal F}_s)$ is the
cokernel of the composition:
$$
\begin{array}{l}
{H^0({\cal F}_1)\otimes \ldots \otimes H^0({\cal F}_s)\lrw
H^0({\cal F}_1\otimes {\cal F}_2)\otimes H^0({\cal F}_3)\otimes \ldots \otimes H^0({\cal F}_s)\lrw }\\
{\lrw H^0({\cal F}_1\otimes \ldots \otimes {\cal F}_s)}\\
\end{array}
$$ \qed

\begin{defin}

Let $V$ be a projective variety and let $|H|$ be a complete unisecant base-point-free linear system
defining a birational map:
$$
\phi_H:V\lrw \ov{V}\subset \P^N
$$
We say that $(V,H)$ is projectively normal or $\Te_V(H)$ normally generated or $\ov{V}$
projectively normal if and only if the natural maps:
$$
\Sym_k(H^0(\Te_V(H)))\lrw H^0(\Te_V(kH))
$$
are surjective for all $k\geq 1$.

\end{defin}

\begin{rem}
{\em Because $H^0(\Te_{P^N}(k))\cong Sim_k(H^0(\Te_{P^N}(1)))$ and $H^0(\Te_{P^N}(1))\cong
H^0(\Te_V(H))$, the following diagram is commutative:
$$
\setlength{\unitlength}{5mm}
\begin{picture}(19,4.5)

\put(0,0){\makebox(0,0){$0$}}
\put(4,0){\makebox(0,0){$H^0(I_{\ov{V},P^N}(k))$}}
\put(10,0){\makebox(0,0){$H^0(\Te_{\P^N}(k))$}}
\put(16.5,0){\makebox(0,0){$H^0(\Te_V(kH))$}}

\put(10,4){\makebox(0,0){$Sim_k(H^0(\Te_V(H)))$}}
\put(0.4,0){\vector(1,0){1.3}}
\put(6.3,0){\vector(1,0){1.3}}
\put(12.5,0){\vector(1,0){1.5}}
\put(10,3.5){\vector(0,-1){3}}
\put(12,3.5){\vector(3,-2){4}}
\put(9.5,2){\makebox(0,0){$\alpha_1$}}
\put(13,0.5){\makebox(0,0){$\alpha_3$}}
\put(15,2.5){\makebox(0,0){$\alpha_2$}}

\end{picture}
$$

Since $\alpha_1$ is surjective, $\coker(\alpha_2)=\coker(\alpha_3)$. Thus, $(V,H)$ is projectively normal
if and only if $s(H,\stackrel{k}{\ldots},H)=0$, for all $k\geq 2$. 

Moreover, we have the following formula to compute the hypersurfaces of degree $k$ containing $\ov{V}$:
$$
h^0(I_{\ov{V},\P^N}(k))=h^0(\Te_{P^N}(k))-h^0(\Te_V(kH))+dim(s(H,\stackrel{k}{\ldots},H))
$$
We will say  that $dim(s(H,\stackrel{k}{\ldots},H))$ is the {\em speciality of $\ov{V}$
respect to hypersurfaces of degree $k$}. 

If the linear system $H$ is very ample,
$V$ and $\ov{V}$ are isomorphic and
$$
H^1(I_{\ov{V},\P^N}(k))\cong s(H,\stackrel{k}{\ldots},H)
$$ }\qed

\end{rem}

\begin{lemma}[Green]\label{green}
Let $\a$,$\b$ be effective divisors on $X$. Let $\b$ be base-point-free. If $h^1(\Te_X(\a-\b))\leq
h^0(\Te_X(\b))-2$ then $s(\a,\b)=0$.
\end{lemma}
{\bf Proof:} Is a particular case of $H^0$-Lemma, \cite{lgreen}. \qed

\begin{lemma}\label{lemma27}
Let $\b_1$,$\b_2$ be divisors on $X$. Let $\a=a_1+\ldots+a_d$ an effective divisor such that
\begin{enumerate}

\item $h^0(\Te_X(\b_1-\a))=h^0(\Te_X(\b_1))-d$,
\item $h^0(\Te_X(\b_2-a_i))=h^0(\Te_X(\b_2))-1$, for $i=1,\ldots,d$.
\item $s(\b_1-\a,\b_2)=0$.
\end{enumerate}
Then, $s(\b_1,\b_2)=0$.
\end{lemma}
{\bf Proof}
See \cite{lange}. \qed

\begin{lemma}\label{penciltrick}
Let $X$ be a smooth curve, ${\cal L}$ an invertible sheaf on $X$ and ${\cal F}$ a torsion-free
$\Te_X-module$. Let $s_1$ and $s_2$ be linearly independent sections of ${\cal L}$ and denote by $V$ the
subspace of $H^0(X,{\cal L})$ they generate. Then the kernel of the cup-product map:
$$
V\otimes H^0(X,{\cal F})\lrw H^0(X,{\cal F}\otimes {\cal L})
$$
is isomorphic to $H^0(X,{\cal F}\otimes {\cal L}^{-1}(B))$ where $B$ is the base locus of the pencil
spanned by $s_1$ and $s_2$
\end{lemma}
{\bf Proof:} See \cite{arbarello}, page 126, Base-point-free Pencil Trick. \qed

\begin{cor}\label{simplepenciltrick}
Let $X$ be a smooth curve. Let $\a$ be a nonspecial divisor such that $|a|$ is a base-point-free
pencil. Then $s(\a,\k)=0$.
\end{cor}
{\bf Proof:} Consider the map:
$$
ker(\alpha)\lrw H^0(\Te_X(\a))\otimes H^0(\Te_X(\K))\lrw H^0(\Te_X(\K+\a))
$$
By applying Lemma \ref{penciltrick} to ${\cal L}=\Te_X(\a)$, $V=H^0(X,{\cal L})$ and ${\cal
F}=\Te_X(\K)$, we obtain that $ker(\alpha)=h^0(\Te_X(\K-\a))=h^1(\Te_X(\a))=0$.
Moreover,
$$
dim(H^0(\Te_X(\a))\otimes H^0(\Te_X(\K)))=2h^0(\Te_X(\K))=2g
$$ 
and
$$
h^0(\Te_X(\K+\a))=a+2g-2-g+1=2g
$$
Then we see that $\alpha$ is surjective (in fact, an isomorphism) and the result follows. \qed

\begin{lemma}\label{importante2}
Let $X$ be a smooth curve of genus $g$. Let $\K$ be the canonical divisor and let $\b$ a
base-point-free nonspecial divisor of degree $b\geq 2g-2$ defining a birational morphism. Then
$s(\b,\K)=0$.
\end{lemma}
{\bf Proof:} Let $\a=a_1+\ldots+a_{b-g-1}$ be a generic effective divisor of degree $b-g-1$, such that 
$|\b-\a|$ is a nonspecial base-point-free linear system. We always can obtain this divisor when $\b$
defines a birational map, by applying the general position theorem.

We will apply Lemma \ref{lemma27}:
$$
\begin{array}{l}
{h^0(\Te_X(\b-\a))=deg(\b-\a)-g+1=h^0(\Te_X(\b))-deg(\a)}\\
{}\\
{h^0(\Te_X(\K-a_i))=h^0(\Te_X(\K))-1\mbox{ because $\K$ is base-point-free}}\\
{}\\
\end{array}
$$
Thus, we only have to prove that $s(\K,\b-\a)=0$. But this is Corollary \ref{simplepenciltrick}.\qed

\begin{lemma}\label{descomponibles}
Let $S=\P(\E_0)$ be a decomposable ruled surface over a smooth curve $X$, with $\E_0\cong \Te_X\oplus
\Te_X(\e)$. Let $|H|=|X_0+\b f|$ be a linear system on $\P(\E_0)$. Then:
$$
\begin{array}{lcccccr}
{s(H,\stackrel{k}{\ldots}
,H)}&{\cong}&{s(\b,\stackrel{k}{\ldots},\b)}&{\oplus}&{s(\b+\e,\b,\stackrel{k-1}{\ldots},\b)}&{\oplus}&{\cdots}\\
{}&{}&{\cdots}&{\oplus}&{s(\b+\e,\stackrel{k}{\ldots},\b+\e)}&{}\\
\end{array}
$$
\end{lemma}
{\bf Proof:} We know that $s(H,\stackrel{k}{\ldots}
,H)$ is the cokernel of the map:
$$
\underbrace{H^0(\Te_S(H)) \otimes \cdots \otimes H^0(\Te_S(H))}_k \lrw
H^0(\Te_S(kH))
$$
Because $\P(\E_0)$ is decomposable, we have the natural isomorphisms:
$$
\begin{array}{l}
{H^0(\Te_S(H))\cong H^0(\Te_X(\b))\oplus H^0(\Te_X(\b+\e))}\\
{}\\
{H^0(\Te_S(kH))\cong \bigoplus_{i=0}^k H^0(\Te_X(k\b+i\e))}\\
\end{array}
$$
We see that $\alpha$ factorizes through the maps
$\alpha_i$:
$$
\setlength{\unitlength}{5mm}
\begin{picture}(20,4.5)

\put(10,4){\makebox(0,0){$\underbrace{H^0(\Te_X(\b))\otimes \ldots \otimes
H^0(\Te_X(\b))}_{k-i}\otimes
\underbrace{H^0(\Te_X(\b+\e)) \otimes \ldots \otimes H^0(\Te_X(\b+\e))}_{i}$}}

\put(10,0){\makebox(0,0){$H^0(\Te_X(k\b+i\e))$}}

\put(10,3.5){\vector(0,-1){3}}
\put(11,2){\makebox(0,0){$\alpha_i$}}

\end{picture}
$$
In this way, $coker(\alpha)=coker(\alpha_1)\oplus \ldots \oplus coker(\alpha_k)$. \qed

We will proof the following theorem:

\begin{teo}\label{principal}
Let $S_{\bb}=\P(\E_{\bb})$ be a canonical ruled surface, such that $H\sim X_0+\b f$ defines a canonical
scroll. Them,

\begin{enumerate}
\item $s(H,H)\cong s(\k,\k)\oplus s(\b,\b)$.

\item $s(H,\stackrel{m}{\ldots},H)\cong
s(\k,\stackrel{m}{\ldots},\k)\oplus s(\b,\stackrel{m}{\ldots},\b)$ for all $m\geq 3$,
except when $g=3$ and $b=4$ or $g=2$ and $b=3$.

\end{enumerate}

\end{teo}
{\bf Proof:} We reduce the proof to check that
$s(\b,\k)=0$, and then we prove this fact in  Lemma
\ref{importante}. 

\begin{enumerate}

\item If $k=2$, by Lemma
\ref{descomponibles} we know that
$$
s(H,H)\cong s(\k,\k)\oplus
s(\k,\b)\oplus s(\b,\b)
$$
so from Lemma \ref{importante} the result follows.

\item If $k>2$, by Lemma
\ref{descomponibles}, we have to proof that
$$
s(\underbrace{\b,\ldots,\b}_i,\underbrace{\K,\ldots,\K}_{m-i})=0
$$
for any $i$, $1\leq i\leq m-1$.

Let us see that $s(\b,\K+\mu \b)=0$, when $\mu\geq
0$. If $\mu=0$, it is Lemma \ref{importante}. If $\mu>0$ we apply Lemma
\ref{green}:
$$
h^1(\Te_X(\K+\mu \b - \b))\leq 1\leq h^0(\Te_X(\b))-2\mbox{ when
$b-g\geq 2$}.
$$

Because we have supposed that $\b$ defines a canonical scroll, the condition 
$b-g\geq 2$ holds except when  $g=2$ and $b=3$ or $g=3$ and $b=4$.

Applying  Lemma \ref{lematecnico} we deduce that:
$$
s(\underbrace{\b,\ldots,\b}_i,\underbrace{\K,\ldots,\K}_{m-i})=s(i\b+\k,\underbrace{\K,\ldots,\K}_{m-i-1})
$$
when $1\leq i \leq m-2$ and
$$
s(\underbrace{\b,\ldots,\b}_{m-1},\k)=s(\b,(m-2)\b+\k)=0
$$

Now, let us see that $s(\k,\lambda \K+\mu \b)=0$, when $\lambda \geq 0$
and
$\mu\geq 1$. We apply Lemma \ref{lemma27}, with $\b_1=\lambda \k+\mu
\b$, $\b_2=\k$ and $\a=\lambda \k+(\mu-1)\b$:
$$
\begin{array}{l}
{h^0(\Te_X(\lambda \K+\mu \b))=h^0(\Te_X(\b))+\deg(\lambda\K+(\mu-1) \b)}\\
{}\\
{h^0(\Te_X(\k-x))=h^0(\Te_X(\k))-1\mbox{ for all $x\in  \lambda\K+(\mu-1)
\b$}}\\
{}\\
{s(\b,\K)=0}\\
\end{array}
$$

Applying Lemma \ref{lematecnico}, we obtain
$$s(i\b+\k,\underbrace{\K,\ldots,\K}_{m-i-1})=s(i\b+(m-i-2)\k,\K)=0$$ when
$1\leq i
\leq m-2$.
\end{enumerate}

\begin{lemma}\label{importante}
Let $X$ be a smooth curve of genus $g$. Let $\K$ be the canonical divisor
and let $\b$ a divisor of degree $b\geq 2g-2$ defining a canonical scroll.
Then $s(\b,\K)=0$.
\end{lemma}
{\bf Proof:} If $X$ is elliptic, $\k\sim 0$ so it is clear that
$s(\b,\k)=0$. Now, we distinguish the hyperelliptic and the nonhyperelliptic cases:

\underline {Case 1: $X$ is hyperelliptic:}

If $|\b|$ defines a birational map it is sufficient to apply Lemma \ref{importante2}.

By Proposition \ref{descripcionb} $|\b|$ does not define a birational map when
$g=2$ and $b=3$ or $g=3$ and
$b=4$. But in these cases $|\b|$ is a base-point-free pencil, so we can apply the Corollary
\ref{simplepenciltrick}. \qed

\underline {Case 2: $X$ is nonhyperelliptic:}

We can apply  Lemma  \ref{importante2} when $|\b|$ defines a birational map. But in this case if
$deg(\b)=2g-2$ $|\b|$ could not verify this condition.

Thus, we will use a different strategy, which will be valid for any divisor $\b$ defining  a canonical
scroll over a nonhyperelliptic curve. We will prove that the following map is a surjection:
$$
H^0(\Te_{S_{\bb}}(H-X_1))\otimes H^0(\Te_{S_{\bb}}(H))\stackrel{\A_0}{\lrw}
H^0(\Te_{S_{\bb}}(2H-X_1))
$$
It holds that
$$
\begin{array}{l}
{H^0(\Te_{S_{\bb}}(H-X_1)))\cong H^0(\Te_{S_{\bb}}(\k f)) \cong H^0(\Te_X(\k))}\\
{H^0(\Te_{S_{\bb}}(H)))\cong H^0(\Te_X(\b))\oplus H^0(\Te_X(\k))}\\
{H^0(\Te_{S_{\bb}}(2H-X_1)))\cong H^0(\Te_X(\b+\k))\oplus H^0(\Te_X(2\k))}\\
\end{array}
$$
so $\A_0$ is a surjection when $s(\b,\k)=0$ and $s(\k,\k)=0$. 

We have the following commutative diagram:
$$
\setlength{\unitlength}{5mm}
\begin{picture}(24,4.5)

\put(8,4){\makebox(0,0){$H^0(\Te_{S_{\bb}}(H-X_1))\otimes H^0(\Te_{S_{\bb}}(H))$}}

\put(8,0){\makebox(0,0){$H^0(\Te_{S_{\bb}}(2H-X_1))$}}

\put(20,2){\makebox(0,0){$H^0(I_{P^{b-g},P^N}(2))$}}

\put(8,3.5){\vector(0,-1){3}}
\put(15,3.5){\vector(2,-1){2}}
\put(17,1.5){\vector(-4,-1){4}}

\put(8.5,2){\makebox(0,0){$\alpha_0$}}
\put(16.1,3.5){\makebox(0,0){$\alpha_1$}}
\put(15.5,0.5){\makebox(0,0){$\alpha_2$}}
\end{picture}
$$ where $\P^N\cong \P(H^0(\Te_{S_{\bb}}(H))^{\vee})$.
Since $\alpha_1$ is a surjection, we only have to proof that $\alpha_2$ is a surjection.
Consider the following exact sequence:
$$
0\lrw \Te_{S_{\bb}}(-C)\lrw \Te_{S_{\bb}}\lrw \Te_C\lrw 0
$$
taking the tensor product with  $\Te_{S_{\bb}}(2H-X_1)$ and applying cohomology we obtain:
$$
\begin{array}{l}
{0\lrw H^0(\Te_{S_{\bb}}(2H-C-X_1))\lrw H^0(\Te_{S_{\bb}}(2H-X_1))\lrw}\\
{\lrw H^0(\Te_C(2H-X_1))\lrw H^1(\Te_{S_{\bb}}(2H-C-X_1))}\\
\end{array}
$$
but $2H-C-X_1$ is a $(-1)$-secant divisor, that is, $\Te_{S_{\bb}}(2H-C-X_1)_x\cong \Te_{P^1}(-1)$, for
all $x\in X$. Because
$H^i(\Te_{P^1}(-1))=0$ for any $i\geq 0$, by the Theorem of Grauert, we deduce that
$H^i(\Te_{S_{\bb}}(2H-C-X_1))=0$ for any $i\geq 0$.
Moreover $\Te_C(2H-X_1)\cong \Te_C(2\k-\Ram)$, so we see that:
$$
H^0(\Te_{S_{\bb}}(2H-X_1))\cong H^0(\Te_C(2\k-\Ram))
$$
In this way to see that $\alpha_1$ is a surjection it is sufficient to check that the following map is a
surjection:
$$
H^0(I_{P^{b-g},P^N}(2))\stackrel{\alpha_3}{\lrw} H^0(\Te_C(2\k-\Ram)
$$
Now, consider the following commutative map:
$$
\setlength{\unitlength}{5mm}
\begin{picture}(24,9)

\put(4,8){\makebox(0,0){$H^0(I_{P^{g-1}\cup C,P^N}(2))$}}

\put(13,8){\makebox(0,0){$H^0(I_{P^{b-g},P^N}(2))$}}

\put(21,8){\makebox(0,0){$H^0(\Te_C(2H-\Ram))$}}

\put(4,4){\makebox(0,0){$H^0(I_{C,P^N}(2))$}}

\put(13,4){\makebox(0,0){$H^0(\Te_{P^N}(2))$}}

\put(21,4){\makebox(0,0){$H^0(\Te_C(2H))$}}

\put(4,0){\makebox(0,0){$H^0(I_{\Ram,P^{b-g}}(2))$}}

\put(13,0){\makebox(0,0){$H^0(\Te_{P^{b-g}}(2))$}}

\put(21,0){\makebox(0,0){$H^0(\Te_{\Ram}(2))$}}

\put(8.3,8.1){\line(0,-1){0.2}}
\put(8.3,4.1){\line(0,-1){0.2}}
\put(8.3,0.1){\line(0,-1){0.2}}

\put(3.9,7.4){\line(1,0){0.2}}
\put(12.9,7.4){\line(1,0){0.2}}
\put(20.9,7.4){\line(1,0){0.2}}

\put(12.9,0.9){\line(1,0){0.2}}
\put(20.9,0.9){\line(1,0){0.2}}

\put(17.5,4.1){\line(0,-1){0.2}}
\put(17.5,0.1){\line(0,-1){0.2}}

\put(8.2,0){\vector(1,0){1.5}}
\put(8.2,4){\vector(1,0){1.5}}
\put(8.2,8){\vector(1,0){1.5}}

\put(16.2,0){\vector(1,0){1.7}}
\put(16.2,4){\vector(1,0){1.7}}
\put(16.2,8){\vector(1,0){1.7}}

\put(4,7.5){\vector(0,-1){3}}
\put(13,7.5){\vector(0,-1){3}}
\put(21,7.5){\vector(0,-1){3}}

\put(4,3.5){\vector(0,-1){3}}
\put(13,3.5){\vector(0,-1){3}}
\put(21,3.5){\vector(0,-1){3}}

\put(4.6,2){\makebox(0,0){$\alpha$}}

\put(17.2,8.5){\makebox(0,0){$\alpha_3$}}

\end{picture}
$$
If we prove that $\alpha$ is a surjection then $\alpha_3$ is a surjective map. We will see this in the
following Lemma: 

\begin{lemma}\label{sobreyectividad}
If $X$ is neither elliptic nor hyperelliptic, then the map:
$$
H^0(I_{C_{\K}}(2))\stackrel{\alpha}{\lrw}
H^0(I_{\Ram,P^{b-g}}(2))
$$ is surjective.
\end{lemma}
{\bf Proof:} We know that $C_{\K}\in \langle 2X_0+\Ram f,2X_1\rangle$.

Let $Q_1$ be a quadric of  $\P^{b-g}$ containing the ramification points  $\Ram$:

\begin{enumerate}

\item If $\ov{X_1}\subset Q_1$, we can take the quadric cone over $Q_1$ of vertex
$\P^{g-1}$,
$Q=\langle \P^{g-1},Q_1\rangle$. This cone
contains the scroll $R$, so it contains $C_{\K}$. Moreover, $Q\cap
\P^{b-g}=Q_1$
so
$\alpha(Q)=Q_1$.

\item If $\ov{X_1}\not\subset Q_1$, $\Ram \subset Q_1$ then $\ov{X_1}\cap
Q_1=\c+\Ram$ where $\c\sim
2\K$. Let
$\c '$ be the set of points corresponding to the divisor $\c$ over $\ov{X_0}$.
Because
$\ov{X_0}$ is not hyperelliptic, it is projectively normal and we can take a quadric
$Q_0$ meeting $\ov{X_0}$ at $\c'$.

Let us take the cones over  $Q_0$ and $Q_1$, and with vertex $\P^{b-g}$ and $\P^{g-1}$
respectively:
$$
\begin{array}{l}
{Cone_0=\langle \P^{g-1},Q_1\rangle}\\
{Cone_1=\langle \P^{b-g},Q_0\rangle}\\
\end{array}
$$
Consider the pencil of quadrics generated by them. Note that any quadric in this pencil contains the
quadric $Q_1$.  Let us study the trace of this pencil on the canonical ruled surface:
$$
\begin{array}{rcl}
{\langle Cone_0,Cone_1 \rangle}&{\lrw}&{H^0(\Te_{S_{\bb}}(2))}\\
{Cone_0}&{\mapsto}&{2X_1+\c f}\\
{Cone_1}&{\mapsto}&{2X_0+\Ram f+\c f}\\
\end{array}
$$
Since $C_{\K}\in \langle 2X_0+\Ram f,2X_1\rangle$, the there exists a quadric
 $Q$ in the pencil containing
 $C_{\K}$ and such that $Q\cap \P^{b-g}=Q_1
$. Therefore $\alpha(Q)=Q_1$. \qed

\end{enumerate}

\begin{rem}
{\em Let us see that in the cases $g=2$, $b=3$ and $g=3$, $b=4$ the Theorem is not true.
Since $s(\b,\k)=0$, we have that $s(\k,\b,\b)=s(\b+\k,\b)$. Furthermore, $h^0(\Te_X(\b))=2$, so we can
compute  $dim(s(\b+\k,\b))$ by using "Base-Point-Free Pencil Trick"
(Lemma \ref{penciltrick}):
$$
dim(s(\b+\k,\b))=h^0(\Te_X(2\b+\k))-2h^0(\Te_X(\b+\k))+h^0(\Te_X(\b+\k-\b))=1
$$
Then we have:
$$
dim(s(H,H,H))=dim(s(\k,\k,\k))+dim(s(\b,\b,\b))+1
$$ }\qed
\end{rem}

Theorem \ref{principal} relates the projective normality of the canonical scroll to the projective
normality  of the directrix curves $\ov{X_0}$ and
$\ov{X_1}$. If a pair $(X,\b)$ is projectively normal and $\b$ is ample, then the linear system
must be very ample (see \cite{mumford}). Then let us see when the divisor $\b$ defining a canonical
scroll is very ample.

\begin{prop}\label{veryample}
Let $\b$ be a nonspecial divisor of degree $b\geq 2g-2$ over a nonhyperelliptic smooth curve
$X$ of genus $g\geq 5$ such that
$\P(\E_{\b})$ is a canonical scroll. Then,

\begin{enumerate}

\item If $b\geq 2g+1$, $|\b|$ is very ample.

\item If $b=2g$ and $g\geq 3$, the generic divisor defining a canonical scroll is very ample.

\item If $b=2g-1$ and $g\geq 5$, the generic divisor defining a canonical scroll is very ample.

\item If $\Cliff(X)\geq 3$, $|\b|$ is very ample.

\end{enumerate}

Furthermore, if $X$ is a generic nonhyperelliptic smooth curve of genus $g>6$, any divisor $\b$
defining a canonical ruled surface is very ample.

\end{prop}
{\bf Proof:} If $b\geq 2g+1$ any divisor of degree $b$ is very ample by Riemann-Roch Theorem.

Let $\b$ be a divisor of degree $b=2g$. Suppose that $|\b|$ is not very ample. Then the
divisor $\b$ must be $\b\sim \k+P+Q$. Since $\b$ defines a canonical ruled surface, $2P+2Q\sim
2\b-2\k$ must be four different points; that is, $2P+2Q\sim P_1+P_2+P_3+P_4$. So there is a
$g^1_4$.

Let $\b$ be a divisor of degree $b=2g-1$. If $|\b|$ is not very ample, $\b-P-Q\sim \K-T$ and
$2\b-2\k\sim 2P+2Q-2T$. Since $2\b-2\k$ must be two different points, $2P+2Q\sim 2T+P_1+P_2$,
so there is a $g^1_4$.

Finally, let $\b$ be a divisor of degree $b=2g-2$. If $|\b|$ is not very ample, $\b-P-Q\sim
\k-T-S$ and $2\b-2\k\sim 2P+2Q-2T-2S$. Since $2\b-2\k$ must be the divisor $0$, $2P+2Q\sim
2T+2S$ and we deduce that there is a $g^1_4$.

Note that a divisor $\b$ which is not very ample is built from $\k$ by using the ramification points of
a $g_4^1$. Then there are a finite number of non very ample divisors $\b$ defining a canonical scroll,
for each $g_4^1$ of $X$.

From this, if the curve $X$ is a generic curve of genus $g>6$ or $\Cliff(X)\geq 3$ then $X$
have not a $g^1_4$, so any divisor $\b$ defining a canonical ruled surface is very ample.

Suppose that $g\geq 5$, by Martens Theorem (\cite{arbarello}, IV, $\S 5$), every component of $W_4^1$
has dimension at most equal to
$1$. Therefore, there is at most a one-dimension family of $g_4^1$. If $b=2g$ (resp.
$b=2g-1$), by Theorem \ref{existenciafuerte}, the family of divisors $\b$ defining a canonical
ruled surface has dimension $4$ (resp. $2$), so the generic divisor is very ample.

Finally, if $g=3$ (resp. $g=4$) and $b=2g$, there is a $3$-dimensional (resp. $4$-dimensional) family
 of divisors defining a canonical scroll. Moreover, there is at most a $2$-dimensional family of
divisors which are not very ample. Therefore, the generic divisor defining a canonical scroll is very
ample.
\qed

\begin{teo}\label{normalidadcanonica}
Let $X$ be a smooth curve of genus $g$. Let $\P(\E_{\b})$ be the
canonical geometrically ruled surface and let $R_{\b}$ be the corresponding canonical scroll.
Then:

\begin{enumerate}

\item If $X$ is not hyperelliptic, then:

\begin{enumerate}

\item If $\deg(\b)\geq 2g+1$, $R_{\b}$ is projectively normal.

\item If $\deg(\b)=2g$, $g\geq 3$, $\b$ is generic then $R_{\b}$ is projectively normal.

\item If $\deg(\b)=2g-1$, $g\geq 5$ and $\b$ is generic, $X$ is not trigonal and $X$ is not of genus $6$
then $R_{\b}$ is projectively normal.

\item If $\Cliff(X)\geq 3$ then $R_{\b}$ is projectively normal.

\end{enumerate}

Furthermore, if $X$ is generic and $g>6$, $R_{\b}$ is projectively normal.

\item If $X$ is hyperelliptic, $R_{\b}$ is not projectively normal.

\item If $X$ is elliptic, $R_{\b}$ is projectively normal.

\end{enumerate}

\end{teo}
{\bf Proof:} We apply Theorem \ref{principal}:

\begin{enumerate}

\item Suppose that $X$ is not hyperelliptic.

Because $\ov{X_0}$ is projectively normal, it is sufficient to check that $\ov{X_1}$ is projectively
normal. 
 
We can suppose that $|\b|$ is very ample in the mentioned cases by Proposition \ref{veryample}.
Our claim follows from Theorem $1$ in \cite{green} and its Corollaries $1.4$ and $1.6$.

\item If $X$ is hyperelliptic. Since $\ov{X_0}$ is not projectively normal, $R_{\b}$ is not
projectively normal.

\item If $X$ is elliptic, $deg(\b)\geq 3$, so $\ov{X_0}$, $\ov{X_1}$ and $R_{\b}$ are projectively
normal.
\qed

\end{enumerate}

\end{document}